\documentclass[twoside,openright,11pt,a4paper,epsf]{article}
\usepackage{latexsym}
\usepackage[utf8]{inputenc}
\usepackage{amsmath}
\usepackage{amsthm}
\usepackage{amsfonts}
\usepackage{amssymb}
\usepackage{epsfig}
\usepackage[all]{xy}
\usepackage{graphicx}

\newcommand{\color}[6]{}
\newcommand{\R}{\mathbb{R}}
\newcommand{\C}{\mathbb{C}}
\newcommand{\D}{\mathbb{D}}
\renewcommand{\P}{\mathbb{P}}

\newcommand{\Z}{\mathbb{Z}}
\newcommand{\Q}{\mathbb{Q}}

\newcommand{\cc}{\mathbb{\mathcal C}}
\newcommand{\cb}{\mathbb{\mathcal B}}
\newcommand{\ca}{\mathcal A}
\newcommand{\ce}{\mathcal E}
\newcommand{\cd}{\mathcal D}
\newcommand{\cv}{\mathcal V}
\newcommand{\cu}{\mathcal U}
\newcommand{\cm}{\mathcal M}

\newcommand{\cl}{\mathcal{L}}
\newcommand{\id}{\text{Id}\,}
\renewcommand{\sp}{\text{Sp}\,}
\newcommand{\nbd}{neighbourhood }
\newcommand{\nbds}{neighbourhoods }
\newcommand{\fonction}[5]
{$$ 
\begin{array}{rcccl}
 #1 & : & #2 & \longrightarrow &#3 \\
    &   & #4 & \longmapsto &#5 
\end{array}
$$}

\newcommand{\rond}[1]{\overset{\vspace*{-1pt}\circ}{#1}}
\newcommand{\priv}{\backslash}
\newcommand{\lra}{\longrightarrow}
\newcommand{\hra}{\hookrightarrow}
\newcommand{\ssi}{\Longleftrightarrow}
\newcommand{\om}{\omega}
\newcommand{\eps}{\varepsilon}
\renewcommand{\phi}{\varphi}
\newcommand{\sdb}{\text{\small SDB}}
\newcommand{\sdbs}{\text{\small SDB}^*}

\newcommand{\msp}{maximal symplectic packing }
\newcommand{\smap}{\,\text{\raisebox{,3cm}{\rotatebox{-90}{$\curvearrowleft$}}}}
\newcommand{\wdt}[1]{\widetilde{#1}}
\newcommand{\cqfd}{\hfill $\square$ \vspace{0.1cm}\\ }
\newcommand{\sbull}{{\tiny $\bullet$ }}
\newcommand{\ds}{\displaystyle}
\newcommand{\im}{\textnormal{Im}\,}

\newtheorem{definition}{Definition}[section]
\newtheorem{thm}{Theorem}
\newtheorem*{thm*}{Theorem}
\newtheorem{prop}[definition]{Proposition}
\newtheorem{lemma}[definition]{Lemma}

\newtheorem{cor}[definition]{Corollary}

\title{ Polarizations and symplectic isotopies.}
\author{Emmanuel Opshtein. \footnote{Partially supported by ANR projects "Floer Power" ANR-08-BLAN-0291-03 and "Symplexe" BLAN06-3-137237}}
\date{}
\begin{document}
\maketitle
\begin{abstract}
The aim of this paper is to explain a link between symplectic isotopies of open objects such as balls and flexibility properties 
of symplectic hypersurfaces. We get connectedness results for spaces of symplectic ellipsoids or  maximal packings  of $\P^2$.
\end{abstract}
\section*{Introduction}
In \cite{biran3}, Biran proved a decomposition theorem for rational k\"ahler  manifolds which has proved useful in many situations 
such as symplectic packings \cite{biran2,moi3}, Lagrangian embeddings \cite{bici1}
  \dots This paper tries to add  symplectic isotopies to this list of applications.
In polarized symplectic manifolds - triples $(M,\om,\Sigma)$  where $\om\in H^2(M,\Z)$ and $\Sigma$ is a symplectic hypersurface 
Poincaré-dual to a (necessarily positive) multiple $k\om$ of the symplectic form -, this decomposition result may be expressed as :
\begin{thm} [Biran]\label{biran} In a polarized closed symplectic manifold $(M,\om,\Sigma)$, there exists a zero-volume closed skeleton in $M$
whose complement is a standard symplectic disc bundle supported by $\Sigma$.
\end{thm} 
 In fact, polarizations of  sufficiently high degree exist on all closed ({\it i.e.} compact, without boundary) rational symplectic manifolds \cite{donaldson}. These manifolds split into a {\it standard symplectic part} - an explicit disc bundle 
 over $\Sigma$ - and a {\it negligible skeleton} (in the K\"ahler case it is even isotropic).
Although no kind of uniqueness can be expected (see \cite{biran3} for a discussion on the skeleton for instance), our 
next theorem explains how to construct many such decompositions - {\it i.e.} with enough freedom - all of whose symplectic parts are
isotopic. In order to state it, let us mention that both $M\priv \Sigma$ and the complement of the zero-section 
$\cl_0$ in the symplectic disc bundle  $\sdb(\Sigma,k)$ are exact symplectic manifolds, so they admit Liouville 
forms. In the statement below,  $\lambda_0$ is a distinguished such  form on $\sdb(\Sigma,k)\priv \cl_0$ (see section \ref{sdbsec}).
Also, by a local embedding of $(X,M)$ into $(Y,N)$ - where $X\subset M$ and $Y\subset N$ - we always mean an embedding of a 
\nbd of $X$ in $M$ into $N$  which induces a diffeomorphism from $X$ to $Y$. 
\begin{thm}\label{flex2}
Let $\phi$ be a local symplectic embedding of $\big(\sdb\,(\Sigma,k),\cl_0\big)$ into $(M,\Sigma)$. Any Liouville form $\beta$ on  $M\priv \Sigma$  which  differs from $\phi_*\lambda_0$ by a smooth $1$-form on $M$ gives rise to a unique symplectic embedding $
\Phi:\sdb(\Sigma,k)\hra ~M$ such that $\Phi_*\lambda_0=\beta$. The embeddings obtained in this way form a contractible space. 
\end{thm}
 The symplectic disc bundles naturally contain ellipsoids or ball packings \cite{moi3}. Theorem \ref{flex2} can be seen as a flexibility statement 
 for embeddings of such objects, which are restricted to a jet condition on a hypersurface. In full generality the flexibility of this condition 
 seems hardly tractable, but in dimension $4$ hypersurfaces are symplectic curves and pseudo-holomorphic techniques come to our 
 rescue. We get for instance :
\begin{thm}\label{isoell}
The space of symplectic embeddings of (closed) balls of fixed radius or ellipsoids $\ce(a_+,a_-)$ with $a_+<\pi$  in $\P^2$
is connected.
\end{thm}
It reproves in particular McDuff's connectedness results on  the space of symplectic  balls of $\P^2$ \cite{mcduff1}, and of 
symplectic ellipsoids \cite{mcduff4} for not too long ellipsoids.  Her approach {\it via} blow-up and inflation is more efficient : 
her results are valid in more general manifolds than only $\P^2$, for more than one ball alone \cite{mcduff2,mcduff3} and 
for any ellipsoid. The present technique may simply be considered as another attempt to approach this still open question 
of symplectic isotopy. In one respect at least it shows interest : compared to the inflation process, the isotopies we produce are more
explicit.  As such, we can  follow them enough  to achieve results for maximal objects. I proved for instance in \cite{moi3}
that all smooth maximal packings of $\P^2$ by two balls - {\it i.e.} smooth embeddings of two closed balls of maximal radii  
($r_1^2+r_2^2=1$) whose images have disjoint interior - have basically the same intersections. In fact, these maximal packings 
are unique up to symplectic isotopies.  
\begin{thm}\label{isomax}
The space of smooth maximal symplectic packings of $\P^2$ by two balls of fixed radii is connected.  
\end{thm}
We do not intend to get the most general results in this paper. For instance theorem \ref{isomax} can be
 easily adapted to five balls and theorem \ref{isoell} to other ellipsoids ({\it e.g.}  $(a_+,a_-)\approx (2\pi,\frac{\pi}{2})$), other spaces such as $S^2\times S^2$ with split symplectic forms (see \cite{moi6} theorem 10) or to results in the directions of Lalonde-Pinsonnault \cite{lapi,lapian}. Our aim is rather 
to present a method, fundamentally based on the fact that since open domains can be made out of symplectic hypersurfaces, these two 
classes of  objects should share some rigidity or flexibility properties. Each  piece of information on the flexibility of some class of 
symplectic hypersurfaces is liable to have an interpretation in terms of packings.

The paper is organized as follows.  We first review Biran's description of the standard disc bundles. We then carry out an alternative proof of theorem \ref{biran} following the lines suggested  in \cite{biran3}. Although this proof gives less information on the skeleton, it is purely symplectic  and it clearly  shows how much flexibility sould be expected for Biran decompositions. It leads to a weak form 
of theorem \ref{flex2},  enough for all subsequent applications.  We come to symplectic isotopies in section \ref{isosec}, 
where we let aside the most delicate point of the method - finding convenient symplectic curves. We deal with this last point in section 
\ref{hardsec} with classical techniques of pseudo-holomorphic curves \cite{gromov,mcsa2}  as well as SFT-like techniques \cite{elgiho,behwz}  similar in spirit to those used by Hind-Kerman \cite{hike}. We finally prove theorem \ref{flex2} in appendix \ref{flex2sec}.

{\small \paragraph{Aknowledgements.} I wish to thank C. Viterbo, F. Bourgeois and A. Oancea  for explaining me how Symplectic Field Theory can be useful  in the kind of problems we are  led to here. This  
paper owes a lot to D. McDuff who  pointed out (and sometimes even 
fixed) inaccuracies in the first versions.}
 
\section{Symplectic disc bundles}\label{sdbsec}
\subsection{Description.}
Let $(\Sigma,\tau)$ be a closed  symplectic manifold with $\tau\in H^2(\Sigma,\Z )$. Let $\pi:\cl \lra \Sigma$ be a line bundle on $\Sigma$ with first chern class $c_1(\cl)=[k\tau]$ endowed with a hermitian metric $g$ and a connection $\nabla$ with curvature 
$R^\nabla=2i\pi k\tau$. Define the transgression $1$-form on the complement of the zero-section $\cl_0$ by :
$$
\begin{array}{l}
\ds \alpha(\frac{\partial}{\partial r})=0, \;\alpha(\frac{\partial}{\partial \theta})=\frac{1}{2\pi},\\
\ds \alpha_{|H^\nabla}=0, 
\end{array}
$$ 
where $r$ is the radial coordinate, $\partial /\partial \theta$ is the infinitesimal generator of the $S^1$-action $e^{i\theta}\cdot$ on $\cl$ and $H^\nabla$ is the horizontal distribution of the connection $\nabla$. With this normalization, a simple computation shows that $d\alpha=-k\pi^*\tau$. Out of this $1$-form, one can make an obviously closed 
$2$-form on $\cl$ :
$$
\om_0:=\pi^*\tau+d(r^2\alpha)=(1-kr^2)\pi^*\tau+dr^2\wedge \alpha.
$$ 
This $2$-form degenerates on the circle bundle $\{r^2=1/k\}$, and therefore restricts to a symplectic form on the $\pi/k$-disc
bundle on $\Sigma$ (that is the set $\{r^2<1/k\}$). Here, $\pi/k$  refers  to  area rather than radius. This disc bundle together with the form $\om_0$ will be called the {\it standard symplectic disc bundle of degree $k$} and denoted by $\sdb(\Sigma,k,\tau)$. 
Of course, this construction relies on a choice of a hermitian metric on $L$, but it can be shown that any such choice 
leads to the same symplectic manifold. This fact justifies the wording "standard" and the omission of the metric on the notation $\sdb(\Sigma,k,\tau)$.  
 Let us point out that these bundles appear naturally in symplectic geometry by Weinstein's theorem \cite{weinstein} :
\begin{thm*}[Weinstein] 
Let $(M,\om)$ be a closed rational symplectic manifold and $\Sigma$ a  symplectic hypersurface of $M$ which is Poincaré-dual to 
an integer $k\om$. Then a \nbd of $\Sigma$ in $M$ is symplectomorphic to a \nbd of $\cl_0$ in $\sdb(\Sigma,k,\om_{|\Sigma})$.
\end{thm*}
When it is clear enough from the context that $\Sigma$ is a hypersurface of a symplectic manifold, we may write $\sdb(\Sigma,k)$ in place 
of $\sdb(\Sigma,k,\om_{|\Sigma}) $.

\subsection{Liouville forms.}  \label{liouvillesec}
The symplectic disc bundles above retract onto their zero-section, so the symplectic form $\om_0$ is exact on the complement of $\cl_0$. Our present aim is to describe some features of their primitives, called Liouville forms, and of their associated vector fields. Throughout this paper, a Liouville form for a symplectic structure $\om$ is a $1$-form $\lambda$ such that $d\lambda=-\om$. Recall that the Liouville vector field associated to  $\lambda$ is defined by $\om(X_\lambda,\cdot)=\lambda$ and is a contracting vector field for $\om$ : $L_{X_\lambda}\om=-\om$, so ${\Phi_{X_\lambda}^t}^*\om=e^{-t}\om$. Let us now fix a symplectic manifold $(\Sigma,\tau)$ and its symplectic disc bundle $\sdb(\Sigma,k,\tau)$ with form $\om_0=(1-kr^2)\pi^*\tau+dr^2\wedge\alpha$. 

The form $\lambda_0:=(1-kr^2)\alpha/k$  is a Liouville form, whose associated vector field is $X_0:=\ds\frac{1-kr^2}{2kr}\frac{\partial}{\partial r}$. This vector field is radial, forward complete (its flow is defined for any positive time), but it explodes at $\cl_0$, where it points outwards. An easy computation shows that if $p\in \sdb(\Sigma,k,\tau)$ is at distance 
$r(p)$ from $\cl_0$, then $\Phi^t_{X_0}(p)$ goes to the boundary of the disc bundle when the time goes forward to $+\infty$
while it reaches $\cl_0$ in finite negative time $\tau(p)=-\ln(1-kr(p)^2)$.
Notice now that if $\vartheta$ is any closed $1$-form on $\Sigma$, the form $\lambda_\vartheta:=\lambda_0+\pi^*\vartheta$ is also a Liouville 
form. Since $\pi^*\vartheta$ vanishes on the fiber, the radial component of the Liouville vector field $X_\vartheta$ associated to 
$\lambda_\vartheta$ is  $X_0$. Moreover, the invariance of $\pi^*\vartheta$ with respect to $r$ ensures that $X_0$ and $X_\vartheta-X_0$ commute. Finally, 
since $\Sigma$ is compact,  $X_\vartheta-X_0$ is complete so the above completeness properties also hold for $X_\vartheta$. In fact - but we will not need this precise statement in the sequel -
$X_\vartheta$ differs from $X_0$ by a horizontal vector field which projects by $\pi$ to the symplectic vector field  $\tau$-dual to $\vartheta$. 
Let us summarize this discussion.
\begin{prop}\label{lvf} Let $(\Sigma,\tau)$ be a closed symplectic manifold with $\tau\in H^2(\Sigma,\Z)$. Any closed 
$1$-form $\vartheta$ on $\Sigma$ gives rise to a Liouville form $\lambda_\vartheta:=\lambda_0+\pi^*\vartheta$ on $\sdb(\Sigma,k)\priv \cl_0$ ($\lambda_0:=(1-kr^2)\alpha/k$) with corresponding vector field $X_\vartheta$ which is forward complete and whose flow satisfies the 
following :
\begin{equation}\label{stop}
\Phi^t_{X_\vartheta}(p)  \lra 
\left\{\begin{array}{l}
 \partial \sdb(\Sigma,k) =\{r^2=1/k\}, \text{ when } t\lra +\infty\\
 \cl_0 =\{r=0\}, \text{ when } t\lra\tau(p):=\ln(1-kr(p)^2).
\end{array} \right.
\end{equation}
\end{prop}
We finish this paragraph with the definition of a convex domain relative to a Liouville form. 
\begin{definition}\label{convexity}
We say that a smoothly bounded compact subset $U\subset \sdb(\Sigma,k)$ is $\lambda_\vartheta$-convex if 
$X_\vartheta$ points outward along $\partial U$.
 \end{definition}
The relevance of this definition lies in that a convex domain $U$ can be recovered from $U\cap \cl_0$ only, by "inflation" along 
$X_\vartheta$. Moreover this inflation process only depends on $X_\vartheta$ inside $U$. For this reason, we call $\cl_0$ a supporting hypersurface
of $U$. 
\subsection{Symplectic disc bundles and ellipsoids.}\label{sdbell}
The role that symplectic disc bundles plays in the applications of this paper comes from their links with ellipsoids. This link was noticed in \cite{moi3}, but we recall it now for self-containedness. Call $\ce(a_1,a_2,\dots,a_n)$
the closed ellipsoid defined in the following way :
$$
\ce(a_1,\dots,a_n)=\left\{ (z_1,\dots,z_n)\in \C^n\;  |\; \pi\frac{|z_1|^2}{a_1}+\dots+\pi\frac{|z_n|^2}{a_n}\leq1 \right\}.
$$
With this definition, the $a_i$ are the areas of the "complex axes" of the ellipsoid, and in symplectic terms they are Eckeland-Hofer capacities of the ellipsoid. A symplectic embedding of an ellipsoid into a symplectic manifold will be simply  called a symplectic ellipsoid. Unless explicitly stated, the ellipsoids are {\it closed} in all the paper. The following proposition shows how they get into the picture (see \cite{moi3} for a proof):
\begin{prop}\label{ellinsdb}
The restriction of a symplectic disc bundle $\sdb(\Sigma,k)$ to a symplectic ellipsoid $\ce(a_1,\dots,a_n)$ of the base is 
an open symplectic ellipsoid $\rond{\ce}(a_1,\dots,a_n,1/k)$. 
\end{prop}
In dimension $4$, that is when $\Sigma$ is of dimension $2$, an ellipsoid of the base is simply an embedded disc in $\Sigma$. 
The proof of proposition \ref{ellinsdb} goes by finding an explicit mapping $\Psi$ between the disc bundle and an ellipsoid. It is easy to prove that smaller {\it straight ellipsoids} (or balls) are identified through $\Psi$ to $\lambda_0$-convex sets. By a   straight  ellipsoid of $\ce(a_1,\dots,a_n)\subset \C^{n}$ we mean the  image of the unit 
ball by a diagonal complex map $(z_1,\dots,z_n)\mapsto (\lambda_1z_1,\dots,\lambda_nz_n)$.
\begin{prop}\label{ballconv}
Let $\ce$ be a symplectic ellipsoid in $\Sigma$ and $\Psi:\ce'\tilde{\lra} \pi^{-1}(\ce)$ be the natural parametrisation of  the restriction of $\sdb(\Sigma,k,\tau)$ to $\ce$ by an ellipsoid.
Then the image of any straight ellipsoid of $\ce'\subset \C^n$ by $\Psi$ is $\lambda_0$-convex as well as $\lambda_\vartheta$-convex provided that the $1$-form $\vartheta$ vanishes  on $\ce$ (this condition is not of cohomological nature).     
\end{prop}

\section{Biran decomposition theorem.}\label{biransec}
\subsection{Proof of theorem \ref{biran}.}
 Let $(M,\om)$ be a closed rational symplectic manifold with 
a symplectic hypersurface $\Sigma$  which is Poincaré-dual to  $k\om$. By Weinstein's \nbd theorem, there is an embedding $\phi$ between open \nbds $\cu,\cv$ of $\cl_0\subset \sdb(\Sigma,k)$ and $\Sigma\subset M$. 
Observe that being Poincaré-dual to $k\om$ implies that $M\priv \Sigma$ is 
an exact symplectic manifold because $\om$ vanishes on any $2$-cycle of $M\priv \Sigma$. The following lemma ensures that one of the $\lambda_\vartheta$ defined in the previous section is compatible with a Liouville form on $M\priv \Sigma$. Precisely,
\begin{lemma}\label{compatibility}
There exists a $1$-form $\vartheta$ on $\Sigma$ such that $\phi_*\lambda_\vartheta$ - defined on $\cv\priv \Sigma$ -
extends to a Liouville form  on $M\priv \Sigma$.
\end{lemma}
\noindent Let us first explain why theorem \ref{biran} is a consequence of this lemma. First, restrict $\cu$ to a 
$\lambda_\vartheta$-convex 
domain if necessary. Let $\beta$ be an extension of $\phi_*\lambda_\vartheta$
to a Liouville form of $M\priv \Sigma$ and $X_\beta$ its dual vector field. Since it coincides with $\phi_*X_\vartheta$ on 
$\cv\priv \Sigma$, $X_\beta$ points outward along $\partial\cv$ so it is forward complete by compactness of $M$, exactly as $X_\vartheta$ is.  Define therefore the map \fonction{\Phi}{\sdb(\Sigma,k)}{M}{x}{\left \{\begin{array}{ll} \phi(x)  & \text{if } x\in \cu, \\ \Phi_{X_\beta}^\tau\circ \phi\circ \Phi_{X_\vartheta}^{-\tau}(x) & \text{if }\Phi_{X_\vartheta}^{-\tau}(x)\in \cu.\end{array}\right.}
Since $\Phi$ transports $X_\vartheta$ to $X_\beta$, it commutes with the flows of these two vector fields inside $\cu$, so there is in fact no need to specify which value of $\tau$ should be considered above : thanks to the convexity of $\cu$, any choice gives the same point  for $\Phi(x)$. Notice 
moreover that the procedure above actually defines $\Phi$ on $\sdb(\Sigma,k)$ because of property (\ref{stop}) : the flow of 
$X_\vartheta$ expands $\cu$ and its image eventually covers the whole of $\sdb(\Sigma,k)$. The map $\Phi$ is obviously a symplectic embedding. 
Finally, notice that $M\priv \im \Phi$ is endowed with a  backward complete (volume-contracting)  Liouville vector field  so its volume must  
vanish by compacity of $M$.\cqfd
Let us finally turn back to the missing lemma.\\
\noindent {\it Proof of lemma \ref{compatibility} :} Pick first a Liouville form $\beta$ on $M\priv \Sigma$. 
If we can find a closed $1$-form $\vartheta$ on $\Sigma$
such that $\beta-\phi_*\lambda_\vartheta$ is exact on $\cv\priv \Sigma$, then the lemma follows. Indeed, $\phi_*\lambda_\vartheta$
then coincides with $\beta+dh$ where $h$ is a smooth function defined on $\cv\priv \Sigma$, and any extension of $h$ 
to $M\priv \Sigma$ provides an extension of $\phi_*\lambda_\vartheta$ to a Liouville form on $M\priv \Sigma$. 

In order to construct the form $\vartheta$, consider a generating family $\{\gamma_0,(\gamma_i)\}$ of the $1$-dimensional homology of $\cv\priv \Sigma$ where $\gamma_0$ is "the small loop around $\Sigma$" parameterized by $\{\phi(e^{i\theta}\cdot p),\; \theta\in[0,2\pi]\}$ and $\gamma_i$  project by $\pi\circ \phi^{-1}$ to simple closed loops $\gamma_i'$ 
 which span $H_1(\Sigma)$.
Let us define $\vartheta$ by requiring that $ [\vartheta]\cdot \gamma_i'=\int_{\gamma_i}(\beta-\phi_*\lambda_0)$. Then $\beta-\phi_*\lambda_\vartheta$ vanishes on each class $[\gamma_i]\in H_1(\cv\priv \Sigma)$ because
$$
\begin{array}{ll}
\ds \int_{\gamma_i}\beta-\phi_*\lambda_\vartheta & = \ds \int_{\gamma_i}\beta - \int_{\gamma_i} \phi_*\lambda_\vartheta \\
 & = \ds \int_{\gamma_i} \beta -\int_{\gamma_i}\phi_*(\lambda_0+\pi^*\vartheta) \\
 & = \ds   \int_{\gamma_i}\big(\beta-\phi_*\lambda_0\big)- \int_{\gamma_i}\phi_*\pi^*\vartheta \\
 & =\ds   \int_{\gamma_i}\big(\beta-\phi_*\lambda_0\big) - \int_{\gamma_i'}\vartheta  = 0 \text{ by construction.}
\end{array}
$$
Observe now that $\gamma_0$ is a torsion element in $H_1(\cv\priv \Sigma)$. Indeed, $[\gamma_0]$ represents the class of the fiber in the $S^1$-bundle $P:=\partial \cv$ over $\Sigma$. Since the Euler class of this bundle is $k[\om]$, the image of $\partial_*$ in the long exact sequence 
$$
\dots\lra H_2(\Sigma) \overset{\partial_*}{\lra} H_1(S^1)\overset{i_*}{\lra} H_1(P)\lra \dots,
$$
  is $[k\om]\big(H_2(\Sigma)\big)\neq \{0\}$. Thus $[\gamma_0]\in H_1(P)=H_1(\cv\priv \Sigma)$ has finite order as claimed and $\beta-\phi_*\lambda_\vartheta[\gamma_0]=0$. Having no period, this one form  is exact in $\cv\priv\Sigma$. \hfill $\square$
\subsection{Comments.}
The proof prompts two remarks. First, we were only able to produce a map which pulls back $\beta$ to some $\lambda_\vartheta$
but not necessarily to $\lambda_0$, which seems in contradiction with theorem \ref{flex2}. The reason for this restriction is 
however purely technical and we get rid of it in the appendix. 
The second remark is that a partial flexibility proceeds already from this  proof of theorem \ref{biran}.
\begin{cor}\label{corflex}
Let $(M,\om,\Sigma)$ be a polarized symplectic manifold, $\phi$ a symplectic embedding of $\sdb(\Sigma,k)$ into $M$ and
  $\vartheta$ a 1-form on $\Sigma$ such that $\phi_*\lambda_\vartheta$ extends to a Liouville form on $M$. 
The set of symplectic embeddings of a $\lambda_\vartheta$-convex domain $U\subset \sdb(\Sigma,k)$ into $M$ which coincide with $\phi$ near $\cl_0$ is contractible. 
\end{cor}
\noindent {\it Proof :} It is enough to prove the statement for those maps which coincide with $\phi$ on a fixed \nbd of $\cl_0$. Then we may as well assume that $U$ contains all $\cl_0$, simply extending all these maps by $\phi$ on some \nbd of $\cl_0$. 
Consider now an element $\psi$ of 
$$
E:=\{\psi:U\lra M\, |\, \psi^*\om=\om_0, \psi=\phi \text{ on a \nbd of } \cl_0\}.
$$
Since $U$ is $\lambda_\vartheta$-convex, it retracts on $\cl_0$ so the cohomological condition on $\vartheta$ implies that $\psi_*\lambda_\vartheta$ 
extends to a Liouville form $\beta$ on $M \priv \Sigma$.  Arguing as in the proof of theorem \ref{biran}, we get a symplectic 
embedding $\Psi_\beta:\sdb(\Sigma,k)\hra M$ with $\Psi_\beta^*\beta=\lambda_\vartheta$. Moreover, because of the $\lambda_\vartheta$-convexity of $U$, the restriction 
of $\Psi_\beta$ to $U$ is $\psi$, whatever extension $\beta$ of $\psi_*\lambda_\vartheta$ we have chosen. This means that 
we have a surjective, obviously continuous map from $E'$ onto $E$, where 
$$
E':=\{\beta\, |\, d\beta=-\om\, ,\, \beta\equiv \phi_*\lambda_\vartheta \text{ on a \nbd of } \cl_0\}.
$$ 
A right-inverse for this map is easily constructed by hand. Since its fibers together with $E'$ are contractible as affine spaces, 
$E$ also is contractible. \hfill $\square$
%\cqfd \vspace{-,5cm}

\section{Applications to isotopy problems in $\P^2$}\label{isosec}
The goal of this section is to exploit corollary \ref{corflex} in some isotopy problems in symplectic geometry. Let us first describe a recipe. Let $(M,\Sigma)$ be a polarized   symplectic manifold $M$ with an embedding of $\sdb(\Sigma,k)$ into $M$ as constructed above.  Let $U$ be a $\lambda_0$-convex domain  of $\sdb(\Sigma,k)$, $\phi:U\hra M$
be the induced embedding and $\psi:U\hra M$ be another symplectic embedding. In order to isotop $\psi$ to $\phi$, proceed as follows : \vspace{,3cm} \\
\sbull\ Find a closed symplectic hypersurface $\Sigma'\subset M$ symplectomorphic to $\Sigma$,
 whose intersection with $\psi(U)$ coincides with $\psi(U\cap \cl_0)$. We say that $\Sigma'$ is a supporting hypersurface of $\psi(U)$. 
 \vspace{,1cm} \\
\sbull\ Isotop $\Sigma'$ to $\phi(\cl_0)$.\vspace{,1cm} \\
{\it Without further information on $\psi$, both previous steps need "hard" techniques, and for this reason the 
present receipe only works in dimension $4$.}\vspace{,1cm} \\
\sbull\ Isotop $\psi$ to $\phi$ on $\Sigma\cap U$, and then even on a \nbd of $\Sigma\cap U$ in $U$. This means finding 
local normal forms, which is often tractable in the symplectic world.\vspace{,1cm} \\
\sbull\  Use finally corollary \ref{corflex} to isotop $\psi$  to $\phi$.\\

Roughly speaking,  the above method reduces the problems of isotopies of open domains from a question on symplectic foliations to a question on one single symplectic hypersurface. Indeed, when isotoping an open domain 
it is {\it a priori} necessary to move full coordinate charts. For instance, Gromov's isotopies between the compactly supported symplectic transformations  of the bidisc  are done through deformations of the images of the grid $\{z_2=c\}$, $\{z_1=c\}$.
McDuff's proofs that the space of embeddings of one ball 
in $\P^2$ is connected consists in straightening the foliation of $\psi(B^4)$ by Hopf discs (the images by $\psi$ of the intersections of $B^4$ with complex lines). 
By contrast, theorem \ref{biran} and its corollary \ref{corflex} show that the data of a single "hard object" - {\it one} symplectic hypersurface - together with  "soft objects" - Liouville forms - provide symplectic coordinates which can already be deformed rather freely. 

In this section, we first focus on the soft part : local isotopies.  We then reprove Gromov's result on the connectedness of the space of symplectic automorphisms of the ball which are the identity near the boundary. Although this result is well-known, we hope that its proof 
will illustrate clearly the above receipe. We finally prove theorem \ref{isoell}  and \ref{isomax}, postponing  however the proof of the needed hard results  to the next section. 
\subsection{Preliminaries : some local normal forms.}
\paragraph{Local ellipsoids.} The first proposition gives local isotopies between symplectic embeddings  of an ellipsoid $\ce:=\ce(1,a)$ which coincide on the Hopf disc $D:=\{z_2=0\}$. A local embedding of $(\ce_a,D)$ is simply an embedding of a \nbd of $D$ in $\ce_a$.
\begin{prop}\label{groloc}
Let $\psi:(\ce_a,D)\to \C^2$ be a local symplectic embedding with $\psi_{|D}=\id$. Then there is a symplectic isotopy $\Phi_t$
which is the identity on $D$ and which isotops $\psi$ to the identity (that is $\Phi_1\circ \psi=\id$). Moreover, the relative version holds : if $\psi_{|\partial \ce_a}=\id$ then one can impose $\Phi_{t|\partial \ce_a}\equiv \id$.
\end{prop}
\noindent {\it Proof :} Since $\psi_{|D}=\id$ the map $\psi$ writes 
$$
\psi(z_1,z_2)=(z_1,A(z_1)z_2) + o(|z_2|),
$$
where  $A(z_1)$ is a symplectic linear map. Consider  a path of maps $A_t:D\to \sp(4)$ such that $A_0=A$ and $A_1=\id$. Impose also that  $A_{t|\partial D}\equiv \id$ when $\psi_{|\partial \ce_a}\equiv \id$, which is possible because $\pi_2(\sp(4))=0$. 
Then there is a path of diffeomorphisms $f_t:(\ce_a,D)\to \R^4$ such that 
\begin{enumerate}
\item $f_0=\phi$, $f_1=\id$,
\item $f_t(z_1,z_2)=(z_1,A_t(z_1)z_2) + o(|z_2|)$,
\item $f_{t|\partial \ce_a}\equiv \id$ if $\psi_{|\partial \ce_a}\equiv \id$.
\end{enumerate}
Since $f_{t*}$ is symplectic on $D$ (and $\partial \ce_a$ in the relative situation), $f_{t*}\om$ differs from $\om$ by a small exact $2$-form near $D$, which can even be chosen to vanish on $D$ (and $\partial \ce_a$). Then Moser's method gives maps $h_t$ which are the identity on $D$ (and $\partial \ce_a$) such that $\phi_t:=h_t\circ f_t$ is a path of symplectic diffeomorphisms. Moreover, 
$h_0=h_1=\id$ because $f_0$ and $f_1$ are symplectic.  \hfill $\square$
\paragraph{The kissing of two spheres along a characteristic.} The next proposition gives a local model, which can be reached 
{\it via} isotopies, for a configuration of two spheres which intersect exactly along one characteristic. It will only be used in the proof of theorem \ref{isomax}. Before stating the result, let 
us describe one such configuration. Consider the two surfaces $S_\pm:=\{|z_1|^2=1\pm |z_2|^2\}\subset \C^2$. Then $S_-$ 
is the $3$-sphere of radius $1$, $S_+\cap S_-=C:=\{(e^{i\theta},0)\}\subset \C^2$ and in fact $S_+$ is also a symplectic sphere locally near $C$.
Consider indeed the complex valued function defined in a \nbd of the unit circle by $\rho(re^{i\theta}):=\sqrt{2-r^2}e^{-i\theta}$. It is an area-preserving involution that preserves the foliation by circles and exchanges the interior and exterior of the unit disc.  Then the 
map  $(\rho(z_1),z_2)$ is a symplectomorphism 
- only defined locally in a \nbd of $C$ -
 which sends $S_-$ to $S_+$.  Let us denote by $\phi_1:=\id$ and $\phi_2:=(\rho(z_1),z_2)$ the local parametrization of $S_-$ and $S_+$ by $S^3\subset \C^2$.  
 \begin{prop}\label{kiss}
Let $(\psi_1,\psi_2)$ be two local symplectic embeddings of $(S^3,C)$ into $(\C^2,C)$ with the properties that 
$\psi_{1|C}=\id$, $\psi_{2|C}(e^{i\theta},0)=(e^{-i\theta},0)$ and $\psi_1(S^3)\cap\psi_2(S^3)=C$. Then there exist symplectic isotopies 
$\Phi_1^t, \Phi_2^t$ with supports in a  \nbd of $C$ such that $\Phi_i^1\circ \psi_i=\phi_i$ near $C$ and 
$\Phi_1^t\circ \psi_1(S^3)\cap\Phi_2^t\circ \psi_2(S^3)=C$ for all time. 
\end{prop} 
\noindent{\it Proof :} Throughout this proof, we will denote by $B_i:=\psi_i(B\cap \cu)$ and $\partial B_i:=\psi_i(S^3)$ where $\cu$ is the \nbd of $C$ on which the $\psi_i$ are defined. 
Consider the orthogonal projections $(D_{\theta})_{\theta \in[0,2\pi[}$ of the discs $\{z_{1}=e^{i\theta},|z_2|<\eps\}$ to the sphere $S^3$. It defines a local foliation $\cd$ of $S^3$ around $C$ by symplectic discs, which 
are tangent to the vertical lines $\{z_{1}=e^{i\theta}\}$. Notice that this foliation is invariant by the diagonal $S^1$-action 
on $\C^2$, whose infinitesimal generator $X$ is tangent to the characteristic foliation of $S^3$. Denote also by $\cd_i=(D_{\theta \,i})$,
$\cd_i'=(D_{\theta\, i}')$, $X_i$ and $X_i'$ the images of  $\cd$ and  $X$ by $\psi_{i}$ and $\phi_{i}$ respectively. \vspace{,1cm}

\noindent {\it Step 1 : Making all the foliations $\cd_i,\cd_i'$ tangent along $C$.} Since all the discs $D_{\theta,i}, \,D_{\theta,i}'$ are transverse to 
$C$, there exists isotopies of foliations $(\cd_i^s)_{s\in [0,1]}$ of $\partial B_i$ whose discs $(D_{\theta,i}^s)$ verify
$$
\left\{
\begin{array}{l}
D_{\theta,i}^0=D_{\theta,i}, \\
D_{\theta,i}(s)\pitchfork C=\{(e^{i\theta},0)\} \hspace{,5cm} \forall s\in[0,1],\\
\cd_i^1 \text{ is tangent to } \cd_i' \text{ along } C.
\end{array}
\right.
$$
Define now an area preserving map $\Phi_i^s:D_{i,0}\to D_{i,0}^s$ and extend it to a self-map of $\partial B_i$ by requiring that $\Phi_i^s$ sends $D_{\theta,i}$ to $D_{\theta,i}^s$ and  preserves the  characteristic foliation  of $\partial B_i$. It is then easy to see that 
$\Phi_i^s$ extends to a symplectomorphism between \nbds of $\partial B_i$. This symplectic isotopy preserves $\partial B_i$
and precisely brings $\cd_i$ to a foliation of $\partial B_i$ which  is tangent to $\cd_i'$. Henceforth, we assume that the four foliations 
$\cd_1,\cd_2,\cd_1',\cd_2'$ are tangent along $C$. \vspace{,1cm}
 
\noindent{\it  Step 2: Making the local {\bf  images} coincide.}   
Thanks to the first step,  all the discs are now tangent so there exists a local diffeomorphism $f$ around $C$ which verifies the following~ : 
\begin{enumerate}
\item $f$ is tangent to the identity on $C$ (that is $f_{*|C}=\id$),
\item $f:D_{0,i}\to D_{0,i}'$ is an area-preserving map,
\item  $f_{*}X_i=X_i'$ (then $f:D_{\theta \, i}\to D_{\theta \, i}'$ is also area-preserving),
\item $f_*$ is symplectic on $\partial B_i$. This is a requirement on the normal derivative of $f$ which is compatible with the condition ($1$),
\end{enumerate}
Now since $f$ is tangent to the identity on $C$, $f_*\om$ differs from $\om$ only by a small exact two form which vanishes 
on $\partial B_i$, provided that we restrict it to a sufficiently small \nbd of  $C$. Moser's path method therefore gives a "correcting map" $h$ which is the identity on $\partial B_1$, $\partial B_2$ and such 
that $h_*f_{*}\om=\om$. The map $\Phi:=h\circ f$ is a symplectic map which transforms $\partial B_i$ to $S_\pm$, so $\Phi\circ \psi_i$ has the same image as $\phi_i$. 

In order to see that $\Phi$ can be obtained by a symplectic isotopy, we use Moser's method once again. Consider a path $(f_t)_{t\in[0,1]}$
of (germs of) diffeomorphisms  near $C$ between $\id$ to $\Phi$ such that all $f_t$ are tangent to the identity on $C$. 
%\begin{enumerate}
%\item[i)] $f_0=\id$, $f_1=\Phi$,
%\item[ii)] $f_{t|C}=\id$,
%\item[iii)] $f_{t*|C}=\id$.
%\end{enumerate}
Again, $f_{t*}\om$ is close to $\om$ near $C$ thanks to condition iii), so we can apply Moser's method and get a path of correction
maps $h_t$ such that $h_t\circ f_t$ is a symplectic isotopy between $\Phi$ and the identity  - recall that $\id$ and $\Phi$ are symplectic 
so $h_0=h_1=\id$.  Notice also that since $h_t$ is a global isotopy, it obviously preserves the feature 
$\partial B_1\cap \partial B_2=C$.  \vspace{.1cm} \\
{\it Step 3 : Reparametrization.} We are now in the situation where $\psi_i(S^3)=\phi_i(S^3)$, $\psi_i(D_\theta)=\phi_i(D_\theta)$
 and $\psi_{i*}X=\phi_{i*}X$. Observe that the
condition $\psi_i(S^3)=\phi_i(S^3)$ together with the fact that these two maps coincide on $C$ impose that 
 their derivatives  only differ by rotations along $C$ :  $\psi_{*i}(p)=\phi_{*i}(p)\circ A^i_p$, where 
$A^i_p(z_1,z_2)=(z_1,e^{i\theta^i_p}z_2)$. Hence rotating the balls allows that $\psi_{*i}(1,0)=\phi_{*i}(1,0)$. Either being more cautious in the previous step, or applying it again at this point of the proof, we can even achieve that $\psi_{i|D_0}=\phi_{i|D_0}$. Since $\Phi^t_X(D_0)=D_t$ and $\psi_{i*}X=\phi_{i*}X$  we thus get $\psi_{i|D_t}=\phi_{i|D_t}$, so $\psi_{i|S^3}=\phi_{i|S^3}$. 
Arguing exactly as for the proof of proposition \ref{groloc}, we get symplectic isotopies from $\psi_i$ to $\phi_i$ which are the 
identity on $\partial B_i$, so they preserve the intersection property $\partial B_1\cap \partial B_2=C$. \hfill $\square$
\subsection{Gromov's theorem}
\begin{figure}[h!]
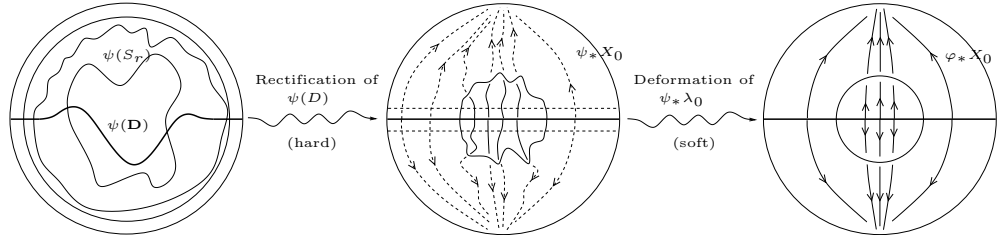

\begin{center} 
\input gromov2.pstex_t
\caption{Proof of Gromov's theorem.}
\label{figgro}
\end{center}\vspace{-.7cm}
\end{figure}

We now prove Gromov's theorem cited above.  Figure \ref{figgro} illustrates the proof.
\begin{thm*}[Gromov] The space of compactly supported symplectic automorphisms of $B^4(1)$ is connected.
\end{thm*}
\noindent {\it Proof :} Let $\psi:B^4(1)\smap$ be a symplectic automorphism which is the identity near $\partial B^4(1)$. Gromov's proof shows 
that the disc $\psi(D):=\psi(\{z_2=0\})$ is symplectic isotopic to $D$ {\it via} a symplectic isotopy with compact support, so we can 
assume that $\psi(D)=D$. Moreover, a Hamiltonian isotopy with support in $D$ transforms $\psi_{|D}$ to the identity (this is an easy $2$-dimensional problem). This isotopy obviously extends to a compactly supported symplectic isotopy, achieving 
$\psi_{|D}=\id$. By proposition \ref{groloc}, we can even assume that $\psi\equiv \id$ on a \nbd of $D$. Since the restriction 
of $\sdb(S^2,1)$ to $S^2\priv\{N\}$ ($N$ is the north pole) is symplectomorphic to $B^4(1)$, we can see $\psi$ as a transformation 
of $\sdb(S^2,1)$ with compact support which is the identity near $\cl_0$ and the fiber over $N$. As in corollary \ref{corflex},
moving $\psi_*\lambda_0$ to $\lambda_0$ {\it via} a path of Liouville forms which fits with $\lambda_0$ wherever $\psi_*\lambda_0=\lambda_0$ gives a path of symplectic automorphisms of $B^4(1)$ which are the identity on $\partial B^4(1)$ and connect 
$\psi$ to the identity.~\cqfd

\subsection{Isotopies of ellipsoids : proof of theorem \ref{isoell}.}\label{isoellsec}
We isotop a symplectic embedding $\psi:\ce(a_+,a_-)\hra \P^2$ with $a_-\leq a_+<\pi$ to a model embedding 
which we describe now. Given a projective line $L$
in $\P^2$, there is an explicit embedding $\Phi$ of $\sdb(S^2,1)$  into $\P^2$, whose complement is 
one point. The restriction of the symplectic disc bundle to the  disc $S^2\priv\{N\}$  is 
symplectomorphic to $B^4(1)$. Our reference embeddings $\phi$ is simply the restriction
of $\Phi\circ s$  to the standard inclusion $\ce(a_+,a_-)\subset B^4(1)$, where $s$ is the coordinate exchange 
$s(z_1,z_2)=(z_2,z_1)$ in $B^4(1)$. Notice that $L$ intersects $\im \phi$ through the small axis $\phi(\{z_1=0\})$ 
of the ellipsoid. Our first task  is to find a 
supporting line for $\psi$ also, that is  a symplectic line of $\P^2$ 
which intersects $\im \psi$ exactly through $\psi(\{z_1=0\})$. 
\begin{lemma} \label{hard}
Every symplectic embedding $\psi:\ce(a_+,a_-)\hra\P^2$ with $a_-\leq a_+<\pi$ 
has a supporting surface of symplectic area $\pi$.  
\end{lemma} 
This lemma is the hardest point of the proof, and it is postponed to the next section. Now since
 any two symplectic lines  are isotopic,  we can assume that the supporting surface is $L$. Moreover, any two discs of the same area in 
 $L$ are isotopic by an ambient Hamiltonian, so we can also assume that $\psi_{|\{z_1=0\}}=\phi_{|\{z_1=0\}}$. Proposition \ref{groloc} actually shows that these embeddings can even fit as germs around $\{z_1=0\}$. By corollary   
 \ref{corflex}, we thus get that $\psi$ is  isotopic to $\phi$.  \hfill $\square$

\subsection{Isotopies of maximal symplectic packings.}
We now explain theorem \ref{isomax}. Let us first recall the setting. A symplectic packing of $\P^2$ by $k$ balls of 
radii $r_1,\dots,r_k$ is a symplectic embedding $\phi:\rond{B}^4(r_1)\sqcup \dots \sqcup \rond{B}^4(r_k)\hra \P^2$. Is is called smooth if 
each $\phi_i:=\phi_{|\rond{B}(r_i)}$ extends to a smooth embedding of the closed ball. It is called maximal if none of the $\phi_i$
can be extended to a larger ball while still defining a packing (that is $\im \phi_i\cap \im \phi_j=\emptyset$). 
Finally, we say that two packings $\phi,\, \psi$ are isotopic if there exist paths of symplectomorphisms $(\Phi_i^t)$
such that 
\begin{itemize}
\item[\sbull] $\Phi_0^i=\id$, $\Phi_i^1\circ \phi_i=\psi_i$
 \item[\sbull] the embeddings $\phi_i^t:=\Phi^t_i\circ \phi_i$ define a packing for each $t\in[0,1]$.
\end{itemize}
When the packing is not maximal and the closure of the balls are disjoint, it is easy to see that this notion of isotopy coincides with 
the notion of ambiant isotopy. However, when the closed balls are not disjoint, this notion allows for more freedom :  breaking intersections or reparametrizing only one of the balls for instance. 
Symplectic packings were first introduced by Gromov, who showed that a symplectic packing of $\P^2$ by two balls 
is subject to a symplectic obstruction : $r_1^2+r_2^2\leq 1$. In \cite{mcpo}, it was proved that this is the only obstruction 
and Karshon constructed an example of smooth maximal $2$-packings. Theorem \ref{isomax} explains that all smooth \msp
 are isotopic to this particular one. We will only address below the problem for two balls of the same capacity $\pi/2$, but the 
 generalization is completely straightforward. 
 
 Since Karshon's example is the packing to which all the others will be isotoped, we describe it briefly now (the description is not Karshon's one). We already pointed out that the restriction of the bundle $\sdb(S^2,1)$ to a disc of area $\pi/2$ is an ellipse $\ce(\pi/2, \pi)$. Covering $S^2$  with two such discs, we get a packing 
 of $\sdb(S^2,1)$ by two ellipsoids $\ce(\pi/2,\pi)$, which contain balls $\cb_1,\cb_2$ of capacity $\pi/ 2$. The embedding $\Phi$ of $\sdb(S^2,1)$ around the line  $L$ already used in the previous section therefore gives a smooth \msp when restricted to $\cb_1,\cb_2$.  Note that the line $L$ is a supporting line of the packing ({\it i.e.} for both balls) and that the two balls intersect along a characteristic circle $C$ of their boundaries (we also write a {\it Hopf circle} in the following).  We denote by $\phi=\{\phi_1,\phi_2\}:=\{\Phi_{|\cb_1},\Phi_{|\cb_2}\}$ this reference maximal packing and by $C_i\subset \partial \cb_i$ the circle $\phi_i^{-1}(C)$.
\begin{figure}[h!]
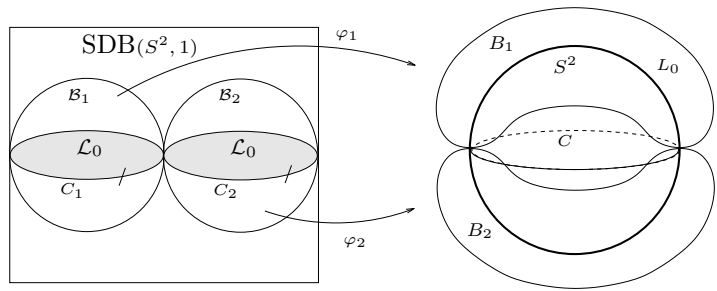

\begin{center} 
\input karshon.pstex_t
\caption{Karshon's maximal packing.}
\label{figkar}
\end{center}\vspace{-.7cm}
\end{figure}

\noindent {\it Proof of theorem \ref{isomax} :} Let $\Psi=\{\psi_1,\psi_2\}$ be a smooth \msp of $\P^2$ by two balls of 
capacity $\pi/ 2$. By \cite{moi3}, we know that $\Psi$ is isotopic to a packing whose balls intersect exactly along one Hopf circle 
of their boundaries. Taking this circle to $C\subset L$ and performing yet another packing isotopy, we can assume that 
$\im \psi_1\cap \im \psi_2=C$, $\psi_i(C_i)=C$ and even that $\psi_{i|C_i}=\phi_{i|C_i}$.  By proposition \ref{kiss}, we can 
even assume that $\psi_i=\phi_i$ near $C$. The intersections $L\cap \psi_i(\cb_i)$ are therefore symplectic discs which coincide
 with the images $\psi_i(\cl_0)$ of the zero section $\cl_0\subset \sdb(S^2,1)$ on their boundary. By Gromov's result,
$\psi_i^{-1}(L\cap \psi_i(\cb_i))$ can be isotoped to $\cl_0$ {\it via} a symplectic isotopy which is the identity on the 
boundary of the ball. So we can safely assume that $\psi_i$ and $\phi_i$ coincide on $\cl_0$. By the relative 
version of proposition \ref{groloc},  we can achieve that $\psi_i=\phi_i$ 
on a \nbd of $\cl_0$. Then $\{\psi_1,\psi_2\}$ can be seen as an embedding $\Psi:\cb_1\cup\cb_2\subset \sdb(S^2,1)\hra \P^2$ which coincides with $\Phi=\{\phi_1,\phi_2\}$ on a \nbd of $\cl_0$. By corollary \ref{corflex}, $\Psi$ and $\Phi$ are 
symplectic isotopic. 
\hfill $\square$
%\cqfd\vspace{-,3cm}

\section{Supporting surfaces of embedded ellipsoids.}\label{hardsec}
We now turn back to the proof of  lemma \ref{hard}. 

\subsection{The case of a ball.} 
The assertion of lemma \ref{hard} when $a_-=a_+=a$ is weaker than what is established in \cite{mcduff1}, namely that {\it all} Hopf discs of 
a ball close to supporting lines. We explain however the proof briefly for the convenience of the reader. It 
 relies on the blow-up procedure introduced by McDuff (see \cite{mcsa}), which we recall now.
Let $\hat\C^2_1:=\{(p,d)\in \C^2\times \P^1,\; p\in d\}$ together with its projections $\pi_1$ on $\C^2$ and $\pi_2$ 
on $\P^1$ (so that $\pi_1:\hat\C^2_1\to \C^2$ is the classical blow-up of $\C^2$ at the origin). Put $E:=\pi_1^{-1}(\{0\})$,
 $L_\delta:=\pi_1^{-1}\big(B^4(\delta)\big)$ and define the symplectic form $\om_\lambda:=\pi_1^*\om_\text{st}+\frac{\lambda}{\pi}\pi_2^*\om_\text{FS}$ on $\hat\C^2_1$. Denote finally by $g:\C^2\priv B^4(\lambda)\to \C^2\priv \{0\}$ the map defined by 
$g(z):=\sqrt{|z|^2-\frac{\lambda}{\pi}}\frac{z}{|z|}$.    Then,
\begin{prop} \label{blowup}
The map $\Phi:=\pi_1^{-1}\circ g$ is a symplectomorphism between $\big(B^4(\lambda+\delta)\priv B^4(\lambda),\om_\text{st}\big)$ and $(L(\delta)\priv L_0,\om_\lambda)$. Moreover, 
\begin{enumerate}
\item[i)] $\Phi^*i$ is $\om_\text{st}$-compatible,
\item[ii)] A disc $D$ intersects $E$ transversally at exactly one point $\xi$ 
if and only if $\Phi^{-1}(D)$ is an annulus one of whose boundary coincides with the circle 
$\Phi^{-1}(\xi)\in \partial B^4(\lambda)$.
\end{enumerate}
\end{prop}
The blow-up procedure can now be summed-up as follows. Given a symplectic embedding $\psi$ of a ball of capacity (slightly larger than) $\lambda$ in $\P^2$, one can use the map $\Phi$ to replace the image of $\psi$ by a \nbd of $L_0$ in $\hat\C^2_1$.
The resulting symplectic manifold is diffeomorphic to the standard complex blow-up of $\P^2$, denoted by $\hat \P^2_1$. Thus the 
embedding $\psi$ defines a symplectic form $\om_\psi$ on $\hat \P^2_1$. It is also possible to extend the standard K\"ahler
structure $i$ defined in a \nbd of $L_0$ to a $\om_\psi$-compatible almost complex structure $J_\psi$ on $\hat \P^2_1$. Notice that the exceptional divisor $E$ is automatically a $J_\psi$-curve. \vspace{,2cm}

\noindent{\it Proof of lemma \ref{hard} ($a_-=a_+=a$):} First blow up  the ball $\psi(B^4(a))$ and consider the induced stuctures 
$\om_\psi,J_\psi$ on $\hat \P^2_1$. In view of the positivity of 
intersection between holomorphic curves and the point ii) above, it is  enough to establish the existence of 
a $J_\psi$-holomorphic curve of $\hat P^2_1$ in the class $L-E$ (where $L$ stands for the projective line and $E$ for the 
exceptional divisor). Indeed, such a curve has area $\pi-a$ and intersects $E$ precisely at one point $\xi$, so its lift by $\Phi$
is a symplectic disc in $\P^2$ with boundary $\partial D=\Phi^{-1}(\xi)$. Gluing the Hopf disc attached to $\Phi^{-1}(\xi)$, we get a topological supporting surface - it may not be smooth along $\Phi^{-1}(\xi)$ - of area $\pi$. Such a topological 
symplectic curve can easily be perturbed to a smooth supporting line for $\psi(B^4(a))$ (\cite{moi3}, lemma 5.2). 

 Let $(\om_t,J_t)_{t\in]0,1]}$ be the symplectic and almost-complex structures on $\hat\P^2_1$ given by the blown-up procedure when carried on for $\psi\big(B^4(ta)\big)$. Denote by $\cm_t$ the moduli space of interest for us :
$$
\cm_t:=\{u:\P^1\lra\hat P^2_1\; |\; du\circ i=J_t\circ du , \; [u]=L-E,\; p\in \im u\},
$$ 
where $p$ is a fixed point of $\hat \P^2_1$ not on the exceptional divisor, and where the moduli-space has to be 
understood modulo reparametrization. Let also be $\cm_t^\eps$ be the same moduli space but for a generic $\eps$-perturbation $J_t^\eps$ of the path $J_t$.
Our task is to understand that $\cm_1$ is not empty. By Alexander's trick, we know that if $\delta$ is small enough, $\psi(B(\delta a))$ is symplectically isotopic to a standard ball in $\P^2$. We can thus assume that $\psi$ and $\phi$ (the standard embedding described in paragraph \ref{isoellsec}) coincide on $B(\delta a)$, which implies that $J_\delta$ is the standard complex structure on $\hat P^2_1$. 
Thus $J_\delta$ is generic for this problem and it is easy to see 
that $\cm_\delta$ consists of exactly one element (the strict transform of a line passing through $p$ and the center of the ball 
in $\P^2$).  If $\cm_1$ is empty, then the space $\{\cm_t^\eps,\; t\in[0,1]\}$ cannot be compact,  so that a bubbling appears by Gromov's  compactness theorem \cite{gromov}. Taking the perturbation $\eps$ small enough, we even see that the bubbling has to appear at some unperturbed almost-complex structure $J_{t_0}$, where $t_0\in ]0,1]$. It means that there is a non trivial decomposition 
$$
L-E=\sum_i A_i\, , 
$$
where $A_i=k_iL+l_i E$ are homology classes of $\hat\P^2_1$ which are represented by  $J_{t_0}$-holomorphic spheres. 
Since holomorphic curves intersect transversally, have positive symplectic area, and noting that $E$ is always $J_t$-holomorphic, we easily see that the decomposition must be of the form 
$$
\begin{array}{ll}
a_0=L-kE& \text{ where } k\geq 2,\\
a_i=k_i E, & \text{ where } k_i\geq 0 \text { and }\sum_{i\geq 1}k_i = k-1.
\end{array}
$$
But we claim that the class $A_0$ is not represented by a $J_{t_0}$-holomorphic curve. Indeed, such a curve would lift to 
a symplectic curve of area $\pi$ in $\P^2$ which coincides with $k\geq 2$ Hopf discs inside the ball $\psi\big(B(t_0 a)\big)$. This symplectic line would therefore have self-intersections which is prohibited by the adjunction inequality. This argument shows that  $\cm_1$ consists of exactly one point. \hfill $\square$
%\cqfd \vspace{-,3cm}

\subsection{Symplectic Field Theory and supporting surfaces.}
This paragraph is aimed at setting the framework for completing the proof of lemma \ref{hard} {\it via} SFT-type arguments. Let 
us first explain the broad idea. Lemma \ref{hard} asserts that one can find a supporting surface of our ellipsoids in the class of a 
symplectic line (of area $\pi$). Following Gromov, we look for this curve among the pseudo-holomorphic spheres for tame almost complex structures. In fact, if we consider a sequence $J_\tau$ of almost complex structures which "stretches the neck" around the boundary of our ellipsoid, the $J_\tau$-holomorphic curves passing through the ellipsoid tend to intersect its boundary along 
closed characteristics when $\tau$  goes to $+\infty$. These curves are therefore good candidates for providing the supporting curves we are looking for. We now review the construction of Symplectic Field Theory  relevant to us and the results we need to prove lemma \ref{hard}.

Let $\psi:\ce(a_+,a_-)\hra \P^2$  be a symplectic embedding of an irrational ellipsoid (that is $a_+,a_-$ are independent 
over $\Q$) and  $\ce$ its image. Then $V:=\partial \ce$ is a hypersurface of contact type in $\P^2$, which has a standard \nbd of the form $(V\times ]-\delta,\delta[,d(e^t\lambda))$ where $\lambda$ is a contact form on $V$ (with associated plane field $\xi$). The Reeb vector field defined on $V$ by $d\lambda(R,\cdot)=0$, $\lambda(R)=1$ and extended by $t$-invariance to $\partial \ce\times ]-\delta,\delta[$ is easily seen to be Hamiltonian. One defines an $\om$-compatible
almost complex structure on $V\times ]-\delta,\delta[$ by 
\begin{equation}\tag{$*$}
\left\{
\begin{array}{l}   
J:\xi\lra \xi \text{ is $\om$-compatible}, \\
\ds J\frac{\partial}{\partial t}=R,
\end{array}
%\tag{\ast}
\right.
\end{equation}
and it can be extended to an $\om$-compatible almost complex structure $J_0$ on $\P^2$. Note that ($*$) can be used to define
an almost-complex structure $\wdt J$ on $V\times \R$. Pulling back $\wdt J$ to $V\times]-\delta,\delta[$ by a map of the form \fonction{\Phi}{V\times]-\delta,\delta[}{V\times]-\tau,\tau[}{(x,t)}{(x,f(t))} where $f'(t)=1$ near $-\delta$ and $\delta$, we get an almost complex structure $J_\tau$ which extends $J_0$ and remains $\om$-compatible. The process of stretching the neck consists in letting $\tau$ go to $\infty$ in this construction. From the complex structure point of view, $\P^2$ then breaks in (at least, see proposition \ref{compactsft}) two pieces bounded by $V$ (the ellipsoid and its complement in our context), both equipped with the almost complex structure $J_\infty$ only defined on the {\it interior} of each piece. What makes this process interesting is a compactness result for sequences of $J_\tau$-holomorphic 
curves due to Bourgeois-Eliashberg-Hofer-Wizocki-Zehnder \cite{behwz}. We sum it up briefly (and not in a strictly rigorous way)
in our context : $V\subset \P^2$ is the boundary of an irrational ellipsoid and the curves are $J_\tau$-holomorphic lines (that is homological to a projective line). Recall that since bubbling cannot appear for this special class of curves, there exists exactly one $J_\tau$-holomorphic line passing through two fixed points for every $\tau<+\infty$. These curves are moreover smoothly embedded spheres.
\begin{prop}\label{compactsft}
Let $\psi:\ce(a_+,a_-)\hra \P^2$ a symplectic embedding of an irrational ellipsoid, $\ce:=\im \psi$, $V:=\partial \ce$,
$\gamma_{\pm}:=\psi(\{z_{2/1}=0\})\cap \partial \ce$ and $J_\tau$ the almost-complex structures defined above. Fix two points in $\P^2$ (at least one lying inside $\ce$), and denote by $S_\tau$ the $J_\tau$-holomorphic line passing through these two points. Then :
\begin{enumerate}
\item[i)] $S_\tau\cap \ce$ converges on any compact set to a collection of $J_\infty$-holomorphic punctured spheres, each of whose ends is asymptotic to some positive multiple of $\gamma_+$ and $\gamma_-$.
\item[ii)] The same happens for $S_\tau\cap {}^c \ce$, except that the ends are asymptotic to negative multiples of $\gamma_+$ and $\gamma_-$. 
\item[iii)] As a whole, in a sense which we do not make precise here (see \cite{behwz}), $S_\tau$ converges to a so-called 
holomorphic building  : the main layers are the ones described in i) and ii) above. The intermediate layers are $\wdt J$-holomorphic punctured spheres on $V\times \R$ asymptotic to positive multiples of $\gamma_1$ and $\gamma_2$ near 
$V\times \{+\infty\}$ and negative multiples near $V\times \{-\infty\}$. All the layers glue together to form a topological sphere. 
\end{enumerate}  
\end{prop}  

\begin{figure}[h!]
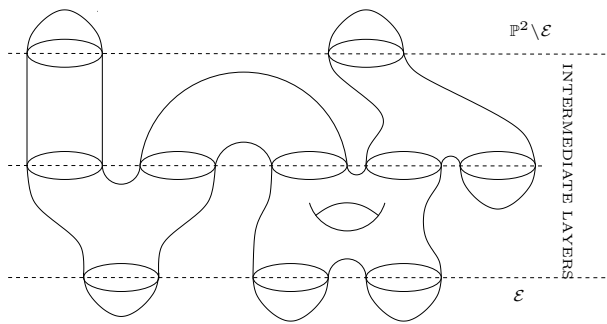

\begin{center} 
\input sft.pstex_t
\caption{A holomorphic building.}
\label{holobuilding}
\end{center}\vspace{-.7cm}
\end{figure}

The previous proposition describes how look like the spheres $S_\tau$ when $\tau$ is large enough. We need to reduce further the possible shapes of these curves. The relevant ingredient is an index formula giving the virtual dimensions 
of the moduli spaces of the involved curves. More specifically, denote by $\cm_\ce(\gamma_1,\dots,\gamma_n,k)$ the space of 
$J_\infty$-spheres with n-punctures in $\ce$, asymptotic to $\gamma_i=l_i \gamma_\pm$ ($l_i>0$) at these punctures, which pass through $k$ fixed points in $\ce$ ($k=1$ or $2$ in our proof), and which are somewhere injective. Similarly, $\cm_{\P^2\priv \ce}((\gamma_i)_{i\leq n},k)$
is the space of somewhere injective $J_\infty$-holomorphic spheres with $n$-punctures in $\P^2\priv \ce$ asymptotic to $\gamma_i=l_i\gamma_\pm$
($k_i<0$) in the homology class of a projective line (relative to $(\gamma_i)_{i\geq n}$) :
$$
\begin{array}{c}
\cm_\ce\big((\gamma_i),k\big):=\left\{
%\begin{array}{ll}
u:S^2\priv\{(x_i)\}\to \ce 
 \left|\begin{array}{l}
 $u$ \text{ is somewhere injective} \\
  du\circ i=J_\infty\circ du \\
 u(\cv(x_i)\priv \{x_i\})\approx \gamma_i, \\   y_1,\dots, y_k\in \im u
 \end{array}\right.\right\}\\
\cm_{\P^2\priv \ce}\big((\gamma_i),k,mL\big):=\left\{u:S^2\priv\{(x_i)\}\to \P^2\priv\ce 
\left|\begin{array}{l}   $u$ \text{ is somewhere injective} \\ du\circ i=J_\infty\circ du\\ u(\cv(x_i)\priv \{x_i\})\approx \gamma_i, \\ y_1,\dots, y_k\in \im u\\ u\in m[L] \text{ rel. } (\gamma_i)_i\} \end{array}\right.\right\}.
\end{array}
$$
We then have (see \cite{elgiho}, proposition 1.7.1 for the dimension of the moduli spaces and \cite{dragnev} for the regularity assertions) :
\begin{prop}\label{dim}
Let $\mu_\text{CZ}(\gamma_i)$ be the Conley-Zehnder index of $\gamma_i\subset \partial \ce\subset \C^2$. Then :
\begin{enumerate}
\item[i)] $\text{VirtDim} \,\cm_\ce\big((\gamma_i),k\big)=-2+\sum (\mu_\text{CZ}(\gamma_i)+1)-2k$,
\item[ii)] $\text{VirtDim}\, \cm_{\P^2\priv \ce}\big((\gamma_i),k,mL\big)=2C_1(mL)-2+\sum (-\mu_\text{CZ}(\gamma_i)+1)-2k=\sum (-\mu_\text{CZ}(\gamma_i)+1)+6m-2-2k$,
\item[iii)] The moduli space $\cm_\ce\big((\gamma_i),k\big)$ or $\cm_{\P^2\priv \ce}\big((\gamma_i),k,mL\big)$ of negative virtual dimension are empty  for a generic choice of $J_0$.
\end{enumerate}
\end{prop}
\subsection{Proof of lemma \ref{hard}.}
Let $\psi:\ce(a_+,a_-)\hra \P^2$ be a symplectic embedding. Denote by $\ce$ its image and $\gamma_+$, $\gamma_-$ 
 the closed characteristics of $\partial \ce$ of $\om$-area $a_+,a_-$ in $\ce$. Let us fix two points $p_1\in\ce$, 
 $p_2\in{}^c\ce$, consider a stretching of the neck $(J_\tau)$  defined as  above, and let $(S_\tau)$
be the path of $J_\tau$-holomorphic lines in $\P^2$ (of symplectic area $\pi$) passing through $(p_1,p_2)$.

When $\tau$ tends to $+\infty$, the spheres $S_\tau$ converge to a $J_\infty$-holomorphic building $S_\infty$ (in the sense of proposition \ref{compactsft}), whose main layers will be denoted in the sequel $S_\infty\cap \ce$ and $S_\infty\cap {}^c\ce$. 
We know that $S_\tau\cap {}^c\ce$ converges to $S_\infty\cap {}^c\ce$ on all compact sets of $\P^2\priv \ce$. 

\begin{lemma} \label{areares}
The surface $S_\infty\cap{}^c\ce$ is made of exactly one connected component.
\end{lemma}
\noindent {\it Proof :} The numbers $n$ of these components is obvisouly  not zero since no compact surface of  $\ce$ is symplectic. Moreover, each connected component of $S_\infty\cap{}^c\ce$
contributes at least  by $\pi$ to the total area of $S_\tau$ for large enough $\tau$ ({\it i.e.} $\ca_\om(S_\tau)\geq n\pi$).
Indeed since $S_\tau$ is close to $S_\infty$ in ${}^c\ce$, it has $n$ connected components, each symplectic. These components  thus have positive symplectic areas on one hand, and are  asymptotic to negative multiples 
$-k_i\gamma_i$ of $\gamma_+,\gamma_-$ on the other hand, so they have area at least $\pi-k_ia_i$. Since the 
complement of these components in $S_\tau$ belongs to a contractible set, the total symplectic area of $S_\tau$ verifies :
$$
\hspace{+2,5cm}\ca_\om(S_\tau) \geq\sum_{i=1}^n (\pi-k_ia_i)+ \sum_{i=1}^n k_ia_i=n\pi. \hspace{3cm} \square
$$   
Call $\gamma_i$ the asymptotes of $S_\infty\cap{}^c\ce$. The curve $S_\infty\cap{}^c\ce$ belongs to the moduli space $\cm_{\P^2\priv\ce}\big((\gamma_i),1\big)$ which must have non-negative virtual dimension by proposition \ref{dim} iv). 
Thus,
$$
\sum (-\mu_\text{CZ}(\gamma_i)+1)+2\geq 0.
$$
Since the Conley-Zehnder indices of the $\gamma_i$ are all no less than $3$ (and $\gamma_-$ is the only orbit with index $3$), we conclude that $S_\infty\cap{}^c\ce$ is a disk asymptotic to $\gamma_-$. This implies that for $\tau$ large enough,
a small Hamiltonian perturbation of $S_\tau\cap{}^c\ce$ makes it a symplectic disk of $\P^2\priv \ce$ 
asymptotic to $\gamma_-$. Closing it by the Hopf disc $\psi(\{z_1=0\})$, we thus get a piecewise smooth symplectic  
curve whose intersection with $\ce$ is exactly $\psi(\{z_1=0\})$. The smoothening of this kind of singular symplectic surfaces is done in \cite{moi3}, lemma 5.2, and it produces a smooth supporting symplectic line for $\ce$.  \hfill $\square$

\subsection{Further supporting polarizations.}
The proof of theorem \ref{isoell} went by finding a supporting line of an ellipsoid passing through his small axis (lemma \ref{hard}). It is however so specific to the situation $a_+<\pi$ that it would be almost dishonest to claim it is really relevant. 
Commonly, the polarizations have high degree, so the ellipsoids which live in the associated disc bundles are very thin and excentric. The needed supporting surface must therefore intersect the ellipsoid we wish to isotop along its big axis instead of 
the  small one. The aim of this paragraph is to show that although it is much harder this way round, finding these 
supporting curves is not always hopeless. We focus on the following :
\begin{thm}\label{isoell2}
The space of symplectic embeddings of ellipsoids $\ce(a_+,a_-)$ in $\P^2$ with $a_+\in]\pi,2\pi[$ and $4a_-<a_+<2\pi<5a_-$ 
is connected.
\end{thm}
Since all symplectic conics are isotopic, the same reasoning as in paragraph  \ref{isoellsec} shows that we only need to find a supporting conic for every embedded ellipsoid $\ce:=\psi\big(\ce(a_+,a_-)\big)$ which intersects $\ce$ along its big axis. 
\begin{lemma}\label{hard2} Every symplectic embedding $\psi:\ce(a_+,a_-)\hra \P^2$  with $a_+\in]\pi,2\pi[$ and 
$\frac{2\pi}{5}<a_-<\frac{a_+}{4}$ admits a supporting conic, that is a symplectic curve of area $2\pi$ whose intersection 
with $\im\psi$ is exactly $\psi(\{z_2=0\})$. 
\end{lemma}
\noindent {\it Proof :} Under the hypothesis of lemma \ref{hard2}, the Conley-Zehnder indices are~:
$$
\begin{array}{c}
\mu_\text{CZ}(\gamma_-)=3,\;\mu_\text{CZ}(2\gamma_-)=5,\;\mu_\text{CZ}(3\gamma_-)=7,\;\mu_\text{CZ}(4\gamma_-)=9,\\\mu_\text{CZ}(\gamma_+)=11,\;\mu_\text{CZ}(5\gamma_-)=13.
\end{array}
$$
Consider a stretching of the neck $(J_\tau)$ around $\partial \ce$. Fix five generic points 
$(p_1,\dots,p_5)$ inside $\ce$. Since bubbling of spheres appear in codimension $2$, a choice of a generic almost-complex structure $J_0$ around $p_1$ insures that for each $\tau$ there is a   $J_\tau$-holomorphic sphere $S_\tau$ in the class of a conic passing through these five points. This sphere is even unique by positivity of intersection. When $\tau$ tends to infinity, $S_\tau$ converges in the sense of \cite{behwz} to a holomorphic building $S_\infty$. Call again $S_\infty\cap\ce$ and $S_\infty\cap{}^c\ce$ the main layers of the building 
respectively inside and outside $\ce$. Reasoning as in lemma \ref{areares} we see that $S_\infty\cap{}^c\ce$ has one 
or two connected components which are punctured spheres. Since $S_\infty$ is a topological sphere, we conclude that 
\begin{itemize}
\item[\sbull] either $S_\infty\cap {}^c\ce$ has two components and the components of $S_\infty\cap \ce$ are disks and 
at most one annulus (the annulus required to connect the two outside connected components may lie in an intermediate 
layer),
\item[\sbull] or $S_\infty\cap{}^c\ce$ has one component and all the components of $S_\infty\cap\ce$ are disks.
\end{itemize}
Using a trick due to Hind-Kerman \cite{hike}, we can get a sharp energy-related restriction on $S_\infty\cap \ce$ : 
\begin{lemma}\label{gamma+} The curve $S_\infty\cap \ce$ is a disk asymptotic to $\gamma_+$.
\end{lemma}
\noindent{\it Proof :} Since $a_-<a_+/5$, the open ellispoid $\rond{\ce}(a_+,a_-)$ contains 
five disjoint closed symplectic balls $\phi_i(B^4(a_+/5-\eps))$, $i=1,\dots,5$.   Starting with an almost complex structure $J_0$
which coincides with ${\psi\circ\phi_i}_*i$ on the images of the $\phi_i$, and $p_1,\dots,p_5$ at the centers of 
these balls, the area of any $J_\infty$-curve passing through $p_1,\dots,p_5$ is no less than $a_+-5\eps$ which is  greater
than $4a_-$ provided that $\eps$ is small enough. The asymptotes of $S_\infty\cap \ce$ must therefore contain at least 
"one $\gamma_+$" or "five $\gamma_-$". Since $\gamma_++\gamma_->5\gamma_->2\pi$, the only possibility is that 
$S_\infty\cap \ce$ has exactly one puncture, asymptotic to $\gamma_+$. \cqfd

Let us finally conclude the proof of lemma \ref{hard2}. Lemma \ref{gamma+} together with the fact 
that  $a_++a_->2\pi$ shows that there can be no intermediate layer. Thus $S_\infty\cap{}^c\ce$ is one disk 
asymptotic to $-\gamma_+$, of area $2\pi-a_+$. Gluing the Hopf disc $\psi(\{z_2=0\})$ to $S_\infty\cap{}^c\ce$, we get 
again a singular symplectic curve, the kind of singularities of which we know how to deal with by \cite{moi3}. We can
smoothen this surface to a symplectic conic coinciding with $\psi(\{z_2=0\})$ inside $\ce$. \hfill $\square$

\appendix
\section{Proof of theorem \ref{flex2}.} \label{flex2sec}
 Theorem \ref{flex2} may be considered as yet another proof 
of theorem \ref{biran}, constructing however the symplectic part of Biran decompositions out of a data with minimal assumptions : 
a polarization $\Sigma$ and a Liouville form on $M\priv \Sigma$  only satisfying a {\it jet condition}. Although this proof is more technical 
and  may not allow for more applications than the proof of section \ref{biransec}, it is more satisfactory in my opinion  since it replaces a 
condition on the germ of a Liouville form - most unnatural in such a flexible world as symplectic geometry - by a much more conceivable 
$1$-jet condition.

 The actual content of theorem \ref{flex2} is contained in the following lemma.
\begin{lemma}\label{lemflex}
Let $\mu$ be a closed form on  $\sdb(\Sigma,k)\priv \cl_0$ which is bounded near $\cl_0$.  Then there is a unique local symplectomorphism $\Psi$ of $(\sdb(\Sigma,k),\cl_0)$ which pulls back $\lambda:=\lambda_0+\mu$ to $\lambda_0$.
\end{lemma}
Indeed, once this fact is established, the map $\Phi:=\phi\circ \Psi$ is a local symplectic embedding of $(\sdb(\Sigma,k),\cl_0)$ into $(M,\Sigma)$ with $\Phi^*\beta=\lambda_0$. The proof of theorem \ref{biran} shows that $\Phi$  extends 
uniquely to a global embedding of $\sdb(\Sigma,k)$ into $M$ with the required property. The last point of the theorem  
concerning the contractibility of this particular kind of embeddings is then obvious since there is a one-to-one correspondence between them and the space of the Liouville forms which are $2$-tangent to $\phi_*\lambda$ - an affine space.

Let us therefore turn back to lemma \ref{lemflex}. Before getting to the core of the proof, let us fix the notation and give the idea.
On each fiber of a the symplectic disc bundle $\sdb(\Sigma,k)$ there are polar coordinates $(r,\theta)$. Unlike the function $r$, the angular
coordinates on the fibers do not fit together to a well-defined global function on of $\sdb(\Sigma,k)\priv \cl_0$ (because of non-vanishing Chern class). However, since nothing in the sequel is affected by this torsion, and in order to avoid irrelevant technicity, we will consider that we have a system of coordinates $(r,\theta,z)$ on $\sdb(\Sigma,k)\priv \cl_0$, where $z$ is a coordinate on $\Sigma$. 
These polar coordinates  naturally extend to coordinates $(r,\theta,z)\in[0,1[\times [0,2\pi[\times \Sigma$ on the 
blow-up of $\cl_0$, which we denote by $\sdbs(\Sigma,k)$. In these coordinates, the flow of the Liouville vector field 
$X_0$ associated to $\lambda_0$ is 
$$
\Phi_0(t,\theta,z):=\Phi_{X_0}^t(0,\theta,z)=(\sqrt{\frac{1-e^{-t}}{k}},\theta,z).
$$ 
It is therefore a diffeomorphism of $\R^*_+\times[0,2\pi[\times\Sigma$ to $\sdb(\Sigma,k)\priv\cl_0$ which extends (only continuously) to a homeomorphism between $\R_+\times[0,2\pi[\times\Sigma$ and $\sdbs(\Sigma,k)$. Since the dual vector field $X$ to $\lambda$ differs from $X_0$ by a bounded vector field, it is tangent to $X_0$ on $\cl_0$ (recall that $X_0$ is of order $r^{-1}$ near $\cl_0$). Although it is not obvious at once, the map $\Phi_X$ defined by $\Phi_X(t,\theta,z):=\Phi_X^t(0,\theta,z)$ also extends to a homeomorphism between $\R_+\times[0,2\pi[\times\Sigma$ and $\sdbs(\Sigma,k)$ which is a diffeomorphism on the complement of $\cl_0$. The diffeomorphism 
$\Psi:=\Phi_X\circ \Phi_0^{-1}$ transports $\lambda_0$ to $\lambda$.
This map $\Psi$ also proceeds from the conjugacy procedure used in the proof of theorem \ref{biran}, but in a more delicate context : if $\tau(p):=-\ln(1-kr^2)$ denotes the time at which the trajectory of $p=(r,\theta,z)$ along $-X_0$
hits $\cl_0$, $\Psi$ is given by : \fonction{\Psi}{\sdbs(\Sigma,k)}{\sdbs(\Sigma,k)}{p=(r,\theta,z)}{\Phi_X^{\tau(p)}(\theta,z)=\Phi_X^{-\ln(1-kr^2)}(\theta,z).}

\noindent {\it Proof of lemma \ref{lemflex} :} Let us first accept the existence and smoothness of the maps $\Phi_X$ and $\Psi$. The derivatives of $\Psi$ in the $(\theta,z)$ directions is clearly the identity on $\cl_0$. Moreover, since $X$ 
is equivalent to $X_{0}=\frac{1-kr^2}{2kr}\frac{\partial}{\partial r}\sim \frac{1}{2kr}\frac{\partial}{\partial r}$ near $\cl_0$, the radial derivative of $\Psi$ is also radial, and the function $f(t):=r(\Phi_X^t(\theta,z))$ satisfy $f'(t)\sim[2kf(t)]^{-1}$
so $f(t)\sim \sqrt\frac{t}{k}$. Therefore, $r(\Psi(r,\theta,z))\sim \sqrt{-\ln(1-kr^2)}\sim \sqrt k r$ so the radial derivative of $\Psi$ is also $1$ and $\Psi$ is tangent to the identity on $\cl_0$.   That $\Psi$ is symplectic is now easy to check. Indeed, $\Psi$ verifies the functional equality 
$\Psi=\Phi^{\tau(p)-\eps}_{X_\lambda}\circ \Psi\circ \Phi^{-\tau(p)+\eps}_{X_0}$, 
so 
$$
\begin{array}{ll}
\Psi^*\om_p & \ds = \Phi^{-\tau(p)+\eps\,*}_{X_0}\Psi^*\Phi^{\tau(p)-\eps\, *}_X\om_p \\
 & \ds = \Phi^{-\tau(p)+\eps\, *}_{X_0}\Psi^*\big(e^{-\tau(p)+\eps}\om_{\Phi^{\tau(p)-\eps}_{X}(p)}\big), \text{ where }  r(\Phi^{\tau(p)-\eps}_{X_\lambda}(p))\underset{\eps\to 0}{\lra} 0 \\
 & = e^{-\tau(p)+\eps} \Phi^{-\tau(p)+\eps\,*}_{X_0}[\om_{\Psi\circ\Phi_{X}^{\tau(p)-\eps}(p})+O(\eps)] =\om_{\Psi(p)}+O(\eps).
\end{array}
$$

Let us turn to the existence and regularity of the involved maps $\Phi_X$ and $\Psi$. In order to investigate the flow of the singular vector field $X$, we compare it with the flow of the {\it regular} vector field $rX$. Since these vector fields are colinear, their flows are related by a 
time rescaling : there exists a function $\tau(t,\theta,z)$ such that 
$$
\Phi_X^t(\theta,z)=\Phi_{rX}^{\tau(t,\theta,z)}(\theta,z).
$$
Differentiating both sides of the previous equality, we get 
$$
\frac{\partial \tau}{\partial t}(t,\theta,z)\, r\big(\Phi^{\tau(t,\theta,z)}_{rX}(\theta,z)\big)=1
$$
or equivalently 
$$
 \int_0^{\tau(t,\theta,z)} r\big(\Phi_{rX}^u(\theta,z)\big)du=t.
$$
\begin{lemma}\label{tfi} There exists a smooth function $f(x,\theta,z)$ on $\R^+\times[0,2\pi[\times \Sigma$  such 
that $\tau(t,\theta,z)=f(\sqrt t,\theta,z)$. Equivalently, 
$$
\int_0^\tau r\big(\Phi_{rX}^u(\theta,z)\big)du=t \ssi \tau=f(\sqrt t,\theta,z).
$$
\end{lemma}

We conclude from this lemma that  $\Phi_X$   extends continuously to $\R^+\times \sdbs(\Sigma,k)$ by the formula $\Phi^t_X(\theta,z)=\Phi_{rX}^{\tau(t,\theta,z)}(\theta,z)=\Phi_{rX}^{f(\sqrt t,\theta,z)}$. Moreover, the map $\Psi$ is  given by 
$$
\Psi(r,\theta,z):=\Phi^{-\ln(1-r^2)}_X(\theta,z)=\Phi_{rX}^{f(\sqrt{-\ln(1-r^2)},\theta,z)}(\theta,z),
$$
so is smooth as a composition of smooth maps (because $\sqrt{-\ln(1-r^2)}$ is indeed smooth on $\R_+$). \cqfd

\noindent {\it Proof of lemma \ref{tfi}:} The vector field $rX=\frac{1-r^2}{2k}\frac{\partial}{\partial r} +rV$ is of class $\cc^1$ on $\sdbs(\Sigma,k)$ and radial on $\cl_0$. Its flow is therefore defined globally on $\sdbs(\Sigma,k)$, and $r\big(\Phi_{rX}^u(\theta,z)\big):\R^+\times [0,2\pi[\times \Sigma\lra [0,1[$ is $\cc^2$-smooth. Moreover, since $rX\sim\frac{1}{2k}\frac{\partial}{\partial r}$ on $\cl_0$, we have $r\big(\Phi_{rX}^u(\theta,z)\big)=\frac{u}{2k}+uh_1(u,\theta,z)$, where $h_1(u,\theta,z)$ is a $\cc^1$-smooth function which vanishes at $u=0$. Now,
$$
\begin{array}{ll}
\ds \int_0^\tau r\big(\Phi_{rX}^u(\theta,z)\big)du& \ds =\int_0^\tau \frac{u}{2k}+uh_1(u,\theta,z)\,du \\
 & \ds =\frac{\tau^2}{4k}+ \tau^2\int_0^1 vh_1(\tau\,v,\theta,z)dv \\
 & \ds = \frac{\tau^2}{4k}\big(1+h_2(\tau,\theta,z) \big)\text{ where }\left\{ \begin{array}{l} h_2 \text{ is } \cc^1-\text{smooth}, \\ h_2(0,\theta,z)\equiv 0 \end{array}\right.\\
 & \ds =\left[\frac{\tau}{2\sqrt k}(1+h_2(\tau,\theta,z))^\frac{1}{2}\right]^2 \\
 & \ds =\left[\tau h_3(\tau,\theta,z)\right]^2,
\end{array}
$$
where $h_3$ is a smooth function which does not vanish at $\tau=0$. Finally, we get that 
$$
\int_0^\tau r\big(\Phi_{rX}^u(\theta,z)\big)du=G(\tau,\theta,z)^2,
$$
where $G$ is a $\cc^1$-smooth function with non-vanishing derivative in the direction of $\tau$ at $\tau=0$. The lemma 
follows from the implicit function theorem. \cqfd

Again this proof has a corollary in terms of embeddings.
\begin{cor}\label{corflex2}
Let $(M,\om,\Sigma)$ be a polarized symplectic manifold,  $\phi$ an embedding of $\sdb(\Sigma,k)$ into $M$ and $\lambda_a$ a Liouville form in $\sdb(\Sigma,k)$ such that $\phi_*\lambda_a$ extends to a Liouville form on $M\priv \Sigma$. Then the space of symplectic embeddings of a $\lambda_a$-convex domain of $\sdb(\Sigma,k)$ whose $1$-jets coincide with the $1$-jet of $\phi$ on $\cl_0$ is contractible. 
\end{cor}

{\footnotesize
\bibliographystyle{abbrv}
\bibliography{bib3.bib}
}

\vspace{2cm}
\noindent Emmanuel Opshtein,\\
Institut de Recherche Mathématique Avancée \\
UMR 7501, Université de Strasbourg et CNRS \\
7 rue René Descartes \\
67000 Strasbourg, France\
opshtein@math.unistra.fr
\end{document}